\UseRawInputEncoding 
\documentclass[reqno,12pt]{amsart}

\usepackage{amsmath}
\usepackage{amssymb}
\usepackage{amsthm}
\usepackage[english]{babel}

\newtheorem{theorem}{Theorem}

\newtheorem{lemma}{Lemma}

\newtheorem{proposition}{Proposition}

\newtheorem{corollary}{Corollary}

\newtheorem{remark}{Remark}

\let\le=\leqslant
\let\ge=\geqslant
\def\bad{\spaceskip=0.33emplus0.6emminus0.15em}

\begin{document}

\title[Concentration Functions]
{Arak Inequalities for Concentration Functions\\
and  the Littlewood--Offord Problem$^{1)}$}


\author[F.~G\"otze]{Friedrich G\"otze$^{*}$}
\author[Yu.S. Eliseeva]{Yulia S. Eliseeva$^{**}$}
\author[A.Yu. Zaitsev]{Andrei Yu. Zaitsev$^{***}$}

\thanks{$^{*}$Fakult\"at f\"ur Mathematik,
Universit\"at Bielefeld, Bielefeld,
Germany; goetze@math.uni-bielefeld.de}
\thanks{$^{**}$St.~Petersburg State University, St.~Petersburg, Russia; pochta106@yandex.ru}
\thanks{$^{***}$St.~Petersburg Department of Steklov Mathematical Institute, St.~Petersburg, Russia
and St.~Petersburg State University, St.~Petersburg, Russia; zaitsev@pdmi.ras.ru}
\thanks{$^{1)}$The paper is supported the SFB 701 in Bielefeld, by Laboratory of
 Chebyshev in St. Petersburg State University (grant of the Government of Russian Federation 11.G34.31.0026),
by grant of St. Petersburg State University
6.38.672.2013, by grants RFBR 13-01-00256, 16-01-00367, by grant NSh-2504.2014.1, and by the Program of Fundamental Researches of Russian
Academy of Sciences "Modern Problems of Fundamental Mathematics".}

\begin{abstract}
Let $X,X_1,\ldots,X_n$ be independent identically distributed random variables.
In this paper we study the behavior of concentration functions of
 weighted sums $\sum_{k=1}^{n}X_ka_k $ depending on the
arithmetic structure of coefficients~$a_k$. The results obtained for the last ten years for the concentration functions of weighted sums play an important role in the study of singular numbers of random matrices. Recently, Tao and Vu proposed
a so-called inverse principle in the Littlewood--Offord problem. We discuss
the relations between this Inverse Principle and a similar principle  for sums of  arbitrarily distributed independent random variables formulated by Arak in the 1980's.
\end{abstract}



\maketitle

\section{Introduction}
\label{s1}

At the beginning of 1980's, Arak~\cite{1}, \cite{2} has published
new bounds for the concentration functions of sums of independent
random variables. These bounds were formulated in terms of the
arithmetic structure of supports of distributions of summands.
Using these results, he has obtained the final solution of an old
problem posed by Kolmogorov~\cite{23}. In this paper, we apply
Arak's results to
 the Littlewood--Offord
problem which was intensively investigated in the last years. We
compare the consequences of Arak's results with recent results of
Nguyen, Tao and Vu~\cite{27}, \cite{28}
and~\cite{35}.

The
 concentration function of a $d$-dimensional random
vector $Y$ with distribution $F=\mathcal L(Y)$ is defined by the
equality
$$
Q(F,\tau)=\sup_{x\in\mathbf{R}^d}\mathbf{P}(Y\in x+ \tau B), \qquad
\tau\ge0,
$$
where $B=\{x\in\mathbf{R}^d\colon \|x\|\le 1/2\}$ is the centered
Euclidean ball of radius~$1/2$.

Let $X,X_1,\dots,X_n$be independent identically distributed
(i.i.d.) random variables. Let $a\,{=}\,(a_1,\dots,a_n)\,{\ne}\,0$, where
$a_k\,{=}\,(a_{k1},\dots,a_{kd})\,{\in}\,\mathbf{R}^d$, $k=1,\dots, n$.
Starting with seminal papers of Littlewood and Offord~\cite{24} and Erd\"os~\cite{13}, the behavior of the concentration
functions of the weighted sums $S_a=\sum_{k=1}^n X_k a_k$.
is studied intensively. In the sequel, let $F_a$ denote the
distribution of the sum $S_a$. The first results were obtained for
the case $\tau=0$ and $d=1$, that is, here the maximal probability $\max_{x\in\mathbf
R}\mathbf{P}\{{S_a=x}\}$. was investigated. For
a detailed history of this part of the problem we refer to a
recent review of Nguyen and Vu~\cite{28}.

In the last ten years, refined concentration results for the
weighted sums $S_a$ play an important role in the study of
singular values of random matrices (see, for instance, Nguyen and Vu~\cite{27}, Rudelson and Vershynin~\cite{31}, \cite{32}, Tao and Vu~\cite{35}, \cite{36}  Vershynin~\cite{38}).

Recently, the authors of the present paper (see~\cite{9},
\cite{10}, and~\cite{12}) improved some of
concentration bounds of the papers~\cite{18}, \cite{31}, \cite{32}, \cite{38}.
These
results reflect the dependence of the bounds on the arithmetic
structure of coefficients~$a_k$ under various conditions on the
vector $a\in {(\mathbf{R}^d)}^n$  
and on the distribution~$\mathcal
L(X)$.

Several years ago, Tao and Vu~\cite{35} (see also~\cite{27}) proposed the so-called inverse principle in
the Littlewood--Offord problem (see~\S\,2). In the
present paper, we discuss the relations between this inverse
principle and similar principles formulated by Arak (see
\cite{1} and~\cite{2})  in his  papers from the 1980's. In the
one-dimensional case, Arak has found a connection of the
concentration function of the sum with the arithmetic structure of
supports of distributions of independent random variables for \textit{arbitrary} distributions of summands.

Apparently the authors of the  publications mentioned above were
not aware of the results from the papers of Arak~\cite{1} and~\cite{2}. Although Arak himself did not use the concept of "inverse principle" in his works, in essence such a principle was there formulated.
It is
related to general bounds for concentration functions of
distributions of sums of independent one-dimensional random
variables. The results were used for the estimation of the rate of
approximation of $n$-fold convolutions of probability
distributions by infinitely divisible ones. Later, the methods
based on Arak's  inverse principle admitted to prove a number of
other important results concerning the rate of infinitely divisible
approximation of convolutions of probability measures.
The problem of estimating this accuracy was formulated by Kolmogorov~\cite{23}. In 1986,
Arak and Zaitsev have published monograph~\cite{3}, containing the above mentioned results, their history and a
discussion of the underlying inverse principle. For the reader's
convenience we include a citation of the relevant passage
concerning this principle from the introduction of monograph~\cite{3}.

{\it ``The concentration functions have turned to be extremely
useful tool in estimating the uniform distance between
convolutions of distributions. They have usually appeared on the
right-hand sides of the corresponding estimates as remainder
terms. However, the general estimates obtained previously for the
concentration functions of $n$-fold convolutions $F^{n}$ were not
sensitive to $Q(F^{n},\tau)$
more rapid than $n^{-1/2}$ in order.

A considerable improvements in the order of estimate
can be achieved by taking into account the structural properties
of the distribution $F$ during the estimation. Already in
considering the example of a distribution~$F$ assigning equal
masses  to points $x_1,\dots, x_m$, it became clear that the
rate of decrease of $Q(F^n,0)$ depends essentially on the mutual
arrangement of these points: $Q(F^{n},0)$ is all the larger, the
more coincidences there are among all possible  numbers of the
form  $\sum_1^m n_k
x_k$, where $n_1,\dots,n_m$ are nonnegative
integers and
$\sum_1^m n_k =n$. Number theory specialists have
known for a long time {\rm (see Freiman~\cite{16})}, that if there are
many such coincidences, then the set
$\{x_1, \dots, x_m\}$ have
an uncomplicated arithmetic structure, in a specific sense.

It turned out that analogous considerations could be
used when the distribution~$F$ is arbitrary, and the
argument~$\tau$ is nonzero: for large $n$ the value of
$Q(F^{n},\tau)$ is essentially greater than zero only if the main
mass of $F$ is concentrated near some finite set~$K$ having a
simple arithmetical structure. It was possible to write this
fairly vague qualitative idea in the form of some new estimates
for concentration functions of distributions of sums of
independent terms.''}

This text is an analogue of descriptions  of the inverse
principles in the papers of Nguyen, Tao and Vu~\cite{27}, \cite{28}
and~\cite{35} (see~\S\,2). A difference being that they restrict
themselves to the classical Littlewood--Offord problem while
discussing the arithmetic structure of the coefficients $a_1,\dots,a_n$ under condition~$Q(F_a, \tau)\ge
n^{-A}$, where $A$ is a positive constant. In this case that one deals with
distributions of sums of non-identically
distributed random vectors of special type only. A further difference is
that, in~\cite{27} and~\cite{28}, the
multivariate case is studied as well.

Nevertheless, there are some \textit{consequences} of Arak's
results which may be interpreted as analogues of
inverse principle for the
Littlewood--Offord problem too. Some of them  have a non-empty
intersection with the results of Nguyen, Tao and Vu~\cite{27}, \cite{28}, \cite{35}, \cite{37} (see
Theorem~\ref{t3}). Moreover, in the monograph~\cite{3}, there
are some structural results (see Theorem~\ref{t4}) implying the assertions which are apparently new in the
Littlewood--Offord problem and have no analogues
in the literature
(see Theorems~\ref{t5} and~\ref{t6}).
We would
like to emphasize that there are of course also some results from~\cite{27}, \cite{28}, \cite{35}, \cite{37} which do not follow from the results of Arak.

Introduce now the necessary notation. The symbol $c$ will be used
for absolute positive constants.
 Note that $c$ can be different in different (or even in the same) formulas.
We will write $A\ll B$, if $A\le c B$. Furthermore, we will use notation $A\asymp B$, if $A\ll B$ and
$B\ll A$. If the corresponding constant depends
on, say, $s$, we write $A\ll_sB$ and $A\asymp_sB$. We denote by
$\widehat F(t)$, $t\in\mathbf R^d$, the characteristic function of
$d$-dimensional distribution~$F$. If $\xi=(\xi_1,\dots,\xi_d)$ is a vector with distribution~$F$, we
denote $F^{(j)}=\mathcal L(\xi_j)$, $j=1,\dots,d$.

For ${x=(x_1,\dots,x_n )\in\mathbf R^n}$, we denote
$$
\|x\|^2=x_1^2+\dots +x_n^2\quad\text{and}\quad |x|= \max_j|x_j|.
$$
Let $E_a$ be the
distribution concentrated at a point $a$. We denote by $[B]_\tau$ the closed $\tau$-neighborhood of
a set $B$ in the sense of the norm $|\,\cdot\,|$.  Products and
powers of measures will be understood in the sense of convolution.
Thus, we write $F^n$ for the $n$-fold convolution of a
measure~$F$. While a distribution~$F$ is infinitely divisible,
$F^\lambda$, $\lambda\ge0$, is the infinitely divisible
distribution with characteristic function $\widehat F^\lambda(t)$.
For a finite set $K$, we denote by $|K|$ the number of
elements~$x\in K$. The symbol $\times$ is used for the
direct product of sets. We write $O(\,\cdot\,)$ if the involved
constants depend on the parameters named "constants" in the
formulations, but not on~$n$.

Let $\widetilde{X}=X_1-X_2$ be the symmetrized random vector,
where $X_1$ and $X_2$ are vectors involved in the definition of
$S_a$ in the Littlewood--Offord problem. In the sequel we use the
notation $G=\mathcal{L}(\widetilde{X})$.

The simplest properties of concentration functions are well
studied (see, for instance,~\cite{3}, \cite{22}, \cite{29}). In particular, it is obvious that
\begin{equation}
\label{eq1}
Q(F,\mu)\le (1+\lfloor\mu/\lambda\rfloor)^d\,Q(F,\lambda),\quad\text{for any }
\mu,\lambda>0,
\end{equation}
where $\lfloor x\rfloor$ is the
largest integer~$k$ that satisfies the inequality $k< x$. Hence,
\begin{equation}
\label{eq2}
Q(F,c\lambda)\asymp_d\,Q(F,\lambda),
\end{equation}
and
\begin{equation}
\label{eq3}
\text{if }\quad Q(F,\lambda)\ll A, \quad \text{then}\quad
Q(F,\mu)\ll A\,(1+\lfloor\mu/\lambda\rfloor)^d.
\end{equation}

Estimating the concentration functions in the
Littlewood--Offord problem, one usually reduces the problem to the
estimation of concentration functions of some symmetric
infinitely divisible distributions. The corresponding statement is
contained in Lemma \ref{l1} below.

For $z\in \mathbf{R}$, introduce the distribution
$H_z$, with the characteristic function
\begin{equation}
\label{eq4}
\widehat{H}_z(t)
=\exp\biggl(-\frac{1}2\sum_{k=1}^n\bigl(1-\cos(\langle t,a_k\rangle z)\bigr)\biggr).
\end{equation}
It is clear that $H_z$ is a symmetric infinitely divisible
distribution.
  Therefore, its characteristic function is positive for all $t\in \mathbf{R}^d$.
For $\delta\ge0$, we denote
\begin{equation}
\label{eq5}
p(\delta)=G\bigl\{\{z\colon |z| > \delta\}\bigr\}.
\end{equation}

\begin{lemma}
\label{l1}
For any $\varkappa,\tau>0$, we have
\begin{equation}
\label{eq6}
Q(F_a, \tau) \ll_d Q(H_1^{p(\tau/\varkappa)}, \varkappa).
\end{equation}
\end{lemma}

According to \eqref{eq3}, Lemma~\ref{l1} implies the following
inequality.

\begin{corollary}
\label{c1}
For any $\varkappa,\tau, \delta>0$, we have
\begin{equation}
\label{eq7}
Q(F_a, \tau) \ll_d(1+\lfloor\varkappa/\delta\rfloor)^dQ(H_1^{p(\tau/\varkappa)}, \delta).
\end{equation}
\end{corollary}

Note that in the case $\delta=\varkappa$ Corollary~\ref{c1}
turns into Lemma~\ref{l1}. Sometimes it is useful to be free in
the choice of $\delta$ in~\eqref{eq7}. In a recent paper of
Eliseeva and Zaitsev \cite{11}, a more general
statement than Lemma \ref{l1} is obtained. It gives useful
bounds if $p(\tau/\varkappa)$ is small, even if
$p(\tau/\varkappa)=0$. The proof of Lemma~\ref{l1} is given
below. It is rather elementary and is based on known properties of
concentration functions. We should note that $H_1^{\lambda}$,
$\lambda\ge0$, is a symmetric infinitely divisible distribution
with the L\'evy
 spectral measure $M_\lambda=(\lambda/4)M^*$, where
 $M^*=\sum_{k=1}^{n}\big(E_{a_k}+E_{-a_k}\big)$.

 Passing in~\eqref{eq6} to the limit, we obtain the following statement (see Zaitsev \cite{48} for details).

\begin{lemma}
\label{l2}
The inequality
\begin{equation}
\label{eq8}
Q(F_a, 0)\ll_d Q(H_1^{p(0)},0)=H_1^{p(0)}\{\{0\}\}
\end{equation}
holds.
\end{lemma}

Lemma~\ref{l1} connects the
Littlewood--Offord problem with general bounds for concentration
functions, in particular with Arak's results. The statement of
Lemma~\ref{l1} is actually the starting point of almost all
recent studies on the Littlewood--Offord problem (usually for
$\tau=\varkappa$, see, for instance,~\cite{18},~\cite{21}, \cite{27}, \cite{31},
\cite{32} and~\cite{38}). More precisely, with the help of Lemma \ref{l3} or its analogues,
the authors of the above-mentioned papers have obtained estimates
of type
\begin{equation}
\label{eq9}
Q(F_a,\tau)\ll_d \sup_{z\ge\tau/\varkappa}\tau^d\int_{|t|\le1/\tau}
\widehat{H}_z^{p(\tau/\varkappa)}(t)\,dt.
\end{equation}
The fact that \eqref{eq1} and~\eqref{eq38} imply that
\begin{align}
&\sup_{z\ge\tau/\varkappa}\tau^d\int_{|t|\le1/\tau}\widehat{H}_z^{p(\tau/\varkappa)}(t)
\,dt\asymp_d\sup_{z\ge\tau/\varkappa}Q(H_z^{p(\tau/\varkappa)},\tau)\nonumber \\
&\qquad=\sup_{z\ge\tau/\varkappa}Q(H_1^{p(\tau/\varkappa)},\tau/z)=
Q(H_1^{p(\tau/\varkappa)},\varkappa),
\label{eq10}
\end{align}
is not used by the authors of these papers. It significantly
hampered the subsequent evaluation of the right-hand side of
inequality~~\eqref{eq9}.

Lemma \ref{l1} reduces the Littlewood--Offord problem to the study of the measure~$M^*$.
 In fact almost all results obtained in solving this problem are formulated in terms of coefficients~$a_j$
 or, equivalently, in terms of properties of the measure~$M^*$.
{This approach does not take into account important information on
 the distribution of the random variable~$X$. In particular, if $\mathcal L(X)$ is standard normal,
 the distribution~$F_a$ is Gaussian with zero mean and covariance operator which is easy to calculate.
 Therefore, there exist bounds for $Q(F_a, \tau)$ which do not follow from any result concerning
 the Littlewood--Offord problem which are discussed in the present paper
 (see, e.g.,~\cite{4}
and~\cite{33}).

In the monograph~\cite{3}, it is also shown
 that if the concentration function of a one-dimensional
  infinitely divisible distribution is large enough, then
 the corresponding L\'evy spectral measure is concentrated
approximately on a set with a special arithmetic structure up to a difference of
 small measure (see Theorems~\ref{t1} and~\ref{t4} below).
Coupled with Lemma~\ref{l1}, these results provide  bounds
in the Littlewood--Offord problem, see Theorems~\ref{t3}, \ref{t5} and~\ref{t6}.

Note that the dependence of the rate of decay of the concentration functions
 of convolutions on the closeness of distributions of summands to some (one-dimensional) lattices  has been
 pointed out even earlier by Mukhin~\cite{26}.
The investigations of  Arak in~\cite{1} and~\cite{2} were motivated  by the ideas
 of Freiman~\cite{16} on the structural theory of set addition. These ideas were used by Nguyen and Vu~\cite{27} and~\cite{28} as well. It should  also be mentioned that Freiman himself
 has used his theory to obtain local limit theorems and bounds for concentration functions
 (see, e.g.,~\cite{8}, \cite{17} and~\cite{25}).

We start now to formulate Theorem~\ref{t1} which is a
one-dimensional Arak type result for infinitely divisible distributions, see~\cite{2}, \cite{3}. Introduce the necessary notations. Let ${\mathbf N}$ be the set
of all positive integers. For any positive integers
$r,m\in{\mathbf N}$ we define $\mathcal{K}_{r,m}$ as the
collection of all sets of the form
\begin{equation}
K=\{\langle{\nu}, h\rangle\colon {\nu}\in \mathbf Z^r\cap V\}\subset\mathbf R,
\label{eq11}
\end{equation}
where $h$ is an arbitrary $r$-dimensional vector, $V$ is an
arbitrary symmetric convex subset of~${\mathbf R}^r$ containing
not more than $m$ points with integer coordinates. That is,
\begin{multline}
\mathcal{K}_{r,m}=\bigl\{\{\langle{\nu},h\rangle\colon {\nu}\in \mathbf Z^r\cap V\}\colon h\in\mathbf R^r,\, V\subset\mathbf R^r,
\\
V=-V,\, V\text{ is convex},\,|\mathbf Z^r\cap V|\le m\bigr\}.
\label{eq12}
\end{multline}

We shall call such sets CGAPs (Convex Generalized Arithmetic
Progressions), by analogy with the notion of GAPs used in the
works of Nguyen, Tao and Vu~\cite{27}, \cite{28} and~\cite{35} (see~\S\,2).

Here, the number~$r$ is the rank and $|\mathbf Z^r\cap V|$ is the volume of a CGAP  in the class   $\mathcal{K}_{r,m}$. It seems natural to call a CGAP from $\mathcal{K}_{r,m}$ proper if all points
$\{\langle{\nu}, h\rangle\colon {\nu}\in \mathbf Z^r\}$
are
disjoint. Notice that, in the definition of the CGAPs, the lattice
${\mathbf Z}^r$ may be replaced by any non-degenerate
$r$-dimensional lattice which may be represented as $\mathbb A\mathbf Z^r$, where $\mathbb
A\colon \mathbf R^r\to \mathbf R^r$ is a non-degenerate linear operator.

For any Borel measure $W$ on ${\mathbf R}$ and $\tau\ge0$ we define
$\beta_{r,m}(W, \tau)$ by the equality
\begin{equation}
\beta_{r,m}(W, \tau)=\inf_{K\in\mathcal{K}_{r,m} }W\{\mathbf R\setminus[K]_\tau\}.
\label{eq13}
\end{equation}

We now introduce  a class  of $d$-dimensional CGAPs $\mathcal{K}_{r,m}^{(d)}$ which consists of all sets of the form $K=\times_{j=1}^d K_j$, where
$K_j\in\mathcal{K}_{r_j,m_j}$, $r=(r_1,\dots,r_d)\in{\mathbf N}^d$,
$m=(m_1, \dots,m_d)\in\mathbf N^d$. We call
$R=r_1+\dots+r_d$ the
rank, and  $|\mathbf Z^{r_1}\cap
V_1|\cdots|\mathbf Z^{r_d}\cap V_d|$ the volume
of~$K$. Here $V_j\subset \mathbf R^{r_j}$ are symmetric convex
subsets from the representation~\eqref{eq11}
for~$K_j$.

The following result is a particular case of Theorem~4.3 of Chapter~II  in~\cite{3}.

\begin{theorem}
\label{t1}
Let $D$ be a one-dimensional
infinitely divisible distribution with characteristic function of
the form $\exp\{\alpha(\widehat
W(t)-1)\}$, $t\in\mathbf R$, where $\alpha>0$ and $W$ is a probability distribution. Let $\tau\ge0$, $r,m\in\mathbf N$. Then 
\begin{equation}
Q(D,\tau)\le c_0^{r+1}\biggl(\frac{1}{m\sqrt{\alpha\beta_{r,m}(W, \tau)}}
+\frac{(r+1)^{5r/2}}{(\alpha\beta_{r,m}(W,\tau))^{(r+1)/2}}\biggr),
\label{eq14}
\end{equation}
where $c_0$ is an absolute constant.
\end{theorem}

Arak~\cite{2} proved an analogue of Theorem~\ref{t1} for sums of i.i.d.\ random variables
(see Theorem~4.2 of Chapter~II  in~\cite{3}). He used this theorem in the proof of the following
remarkable result:
\textit{There exists a universal constant~$C$ such that for
any one-dimensional probability distribution~$F$ and for any
positive integer~$n$ there exists an infinitely divisible
distribution~$D_n$ such that
\begin{equation}
\label{eq15}
\rho(F^n, D_n)\le C\,n^{-2/3},
\end{equation}
where $\rho(\,\cdot\,,\,\cdot\,)$ is the classical
Kolmogorov uniform distance between corresponding distribution
functions.}

This gave the final solution to the long-standing problem stated by
Kolmogorov~\cite{23} in the 1950's (see~\cite{3} for
the history of this problem). Note that the rate of
approximation in~\eqref{eq15} is much better than the rate of
approximation in the well-known
 Berry--Ess\'een theorem. Moreover, the distribution~$F$ is
{\it arbitrary}, no moment type assumptions are imposed. In
addition, this result is in a natural sense unimprovable
(see~\cite[Chapter~VIII]{3}).

Below we will use  the condition
\begin{equation}
G\{\{x\in\mathbf{R}\colon C_1<|x| < C_2\}\}\ge C_3,
\label{eq16}
\end{equation}
where the values of $C_1, C_2, C_3$ will be specified in the
formulations below. Lemma~\ref{l1} and Theorem~\ref{t1} imply the
following Theorem~\ref{t2}.

\begin{theorem}
\label{t2}
Let $\varkappa, \delta>0$,  $\tau\ge0$, and let $X$ be a real random variable satisfying
 condition~\eqref{eq16} with~$C_1=\tau/\varkappa$, $C_2=\infty$
and $C_3=p(\tau/\varkappa)>0$. Let $d=1$, $r,m\in\mathbf N$. Then
\begin{equation}
Q(F_a, \tau)\le c_1^{r+1}(1+\lfloor\varkappa/\delta\rfloor)\biggl(\frac{1}{m\sqrt{\beta_{r,m}(M_0, \delta)}}
+\frac{(r+1)^{5r/2}}{(\beta_{r,m}(M_0, \delta))^{(r+1)/2}}\biggr),
\label{eq17}
\end{equation}
and, for $\tau=0$,
\begin{equation}
Q(F_a, 0)\le c_1^{r+1}\,\biggl(\frac{1}{m\sqrt{\beta_{r,m}(M_0, 0)}}
+\frac{(r+1)^{5r/2}}{(\beta_{r,m}(M_0, 0))^{(r+1)/2}}\biggr),
\label{eq18}
\end{equation}
where $M_0=\frac{p(\tau/\varkappa)}4M^*$, $M^*=\sum_{k=1}^n(E_{a_k}+E_{-a_k})$ and $c_1$ is an absolute constant.
\end{theorem}

In order to prove Theorem~\ref{t2}, it suffices to apply
Corollary~\ref{c1}, Lemma~\ref{l2} and Theorem~\ref{t1}, and to note that $H_1^{p(\tau/\varkappa)}$ is an infinitely
divisible distribution with L\'evy
 spectral measure $M_0$.
 Introduce also
$M=\sum_{k=1}^nE_{a_k}$.  It is obvious that $M\le M^*$ and
$\beta_{r,m}(M, \delta)\le\beta_{r,m}(M^*, \delta)$.

Theorem \ref{t3} follows from Theorem~\ref{t2}. The
conditions of this theorem are weaker than those
used in the results of Nguyen, Tao and Vu~\cite{27}, \cite{28} and~\cite{35}.  In \S\,{2},
we compare Theorem~\ref{t3} with these results.

\begin{theorem}
\label{t3}
Let $d\ge1$, $0 < \varepsilon \le 1$, $0 < \theta \le 1$,
$A>0$, $B>0$, $C_3>0$ be constants and   $\tau=\tau_{n} \ge 0$
be a
parameter that may depend on $n$. Let $X$ be a real random variable satisfying
 condition~\eqref{eq16} with~$C_1=1$, $C_2=\infty$ and $C_3\le p(1)$.
Suppose that
$a=(a_1,\dots,a_n)\in {(\mathbf{R}^d)}^n$ is a multivector in~$\mathbf R^d$ such
that $q_j=Q(F_a^{(j)},\tau)\ge n^{-A}$,
$j=1,\dots,d$, where $F_a^{(j)}$ are distributions of coordinates of the vector~$S_a$.
Let $\rho_n$ denote  a non-random sequence satisfying $n^{-B}\le\rho_n\le1$.
Then, for any number $n'\in{\mathbf N}$ between $\varepsilon n^\theta$ and $n$,
 there exists a CGAP $K$ such that

{\rm 1)}~At least $n-dn'$ elements of\/ $a$ are
$\tau\rho_n$-close to $K$ in the norm $|\,\cdot\,|$ $($this means that, for these elements~$a_j$, there exist $y_{j}\in K$ such
that $|a_{j}-y_{j}|\le\tau\rho_n)$;

{\rm 2)}~$K$ has small rank $R=O(1)$, and small cardinality
\begin{equation}
|K|\le\prod_{j=1}^d\max\bigl\{O\bigl(q_j^{-1}\rho_n^{-1}(n')^{-1/2})\bigr),
1\bigr\}.
\label{eq19}
\end{equation}
\end{theorem}

\begin{remark}
\label{r1}\rm
In Theorem~\ref{t3}, the CGAP $K$ may be non-proper.
\end{remark}

Theorem~\ref{t1} has been proved for one-dimensional situations
and thus initially allows  us
 to  prove Theorem~\ref{t3} for $d=1$ only.
  However, we will show that this one-dimensional version of
Theorem~\ref{t3} provides sufficiently rich  arithmetic
properties for the set $a=(a_1,\dots,a_n)\in{(\mathbf{R}^d)}^n$
in the multivariate case as well. To this end it suffices to apply the
one-dimensional version of Theorem~\ref{t3} to the
distributions
$F_a^{(j)}$, $j=1,\dots,d$. Notice that the
condition $Q(F_a,\tau)\,{\ge}\,n^{-A}$ 
implies that $Q(F_a^{(j)}, \tau)\,{\ge}\allowbreak n^{-A}$, 
$j=1,\dots,d$, since $Q(F_a^{(j)}, \tau)\ge Q(F_a, \tau)$.

Theorem \ref{t2} has non-asymptotic character, it is more
general than Theorem~\ref{t3} and gives information about the
arithmetic structure of  $a=(a_1,\dots,a_n)$ without assumptions
like
$q_j=Q(F_a^{(j)},\tau)\ge n^{-A}$, $j=1,\dots,d$. Notice that in the
asymptotic Theorems~\ref{t3}, \ref{t12} and~\ref{t13}, where $n\to\infty$,
the elements~$a_{j}$ of the multivector~$a$ may depend on~$n$.

Below we formulate another one-dimensional result of Arak (see
Theorem~\ref{t4}). Theorem~\ref{t4} will allow us to prove
another inverse principle type result in the
Littlewood--Offord problem.

For any $r\in\mathbf N$ and $u=(u_1,\dots, u_r)\in {({\mathbf
R}^d)}^r$, $u_j\in \mathbf R^d$, $j=1,\dots,r$, we introduce
the set
\begin{equation}
{K}_1(u)=\biggl\{\sum_{j=1}^r n_j u_j\colon n_j\in \{-1,0,1\} \text{ ¤«п } j=1,\dots,r\biggr\}.
\label{eq20}
\end{equation}
 Define also collection of sets
\begin{equation}
\mathcal{K}_{r}^{(d)}=\bigl\{{K}_1(u)\colon u=(u_1,\dots, u_r)\in
{(\mathbf R^d)}^r\bigr\}.
\label{eq21}
\end{equation}
It is easy to see that the set $K_1(u)$ is symmetric GAP
of rank~$r$ and volume~$3^r$ (see~\S\,\ref{s2}).

The following Theorem~\ref{t4} is Theorem~3.3 of Chapter~II of the monograph~\cite{3}.
 It follows directly from the results of
Arak~\cite{1}.

\begin{theorem}
\label{t4}
{\bad Let $D$ be a one-dimensional
infinitely divisible distribution with characteristic function of
the form $\exp\{\alpha(\widehat
W(t)-1)\}$, $t\in\mathbf R$, where $\alpha>0$,  and~$W$ is a  one-dimensional probability
distribution. Let $\tau\ge0$ and $\gamma=Q(D,
\tau)$}. 
Then there
exist $r\in\mathbf N$ and numbers $u_1,\dots,
u_r\in {\mathbf R}$ such that
\begin{equation}
r\ll|\ln\gamma|+1
\label{eq22}
\end{equation}
\vskip-2mm 

\noindent and
\begin{equation}
\alpha W\{\mathbf R^d\setminus[K_1(u)]_\tau\}\ll(|\ln \gamma|+1)^3,
\label{eq23}
\end{equation}
\vskip-2mm 

\noindent where $u=(u_1,\dots, u_r)\in \mathbf R^r$.
\end{theorem}

 Theorem~\ref{t4} was also used for estimation of the rate of
infinitely divisible approximation of convolutions of probability
distributions (see~\cite{1}, \cite{3},
\cite{5}--\cite{7}, \cite{39}--\cite{49}).

In particular, Zaitsev (see~\cite{49}) solved another problem
considered in the 1950s by  Kol\-mo\-go\-rov~~\cite{23}.
He managed to
get the correct order of the accuracy of infinitely divisible
approximation of distributions of sums of independent random
variables, the distribution of which are concentrated on the short
intervals of length $\tau\le1/2 $ to within a small probability~$p$.
It was found that the accuracy of approximation in the L\'evy
metric has order
$ p+\tau\ln(1/\tau)$, which is much
more accurate than the initial result of Kolmogorov
$p^{1/5}+\tau^{1/2}\ln(1/\tau)$, and
later obtained results of other authors. As approximation, the
so-called accompanying infinitely divisible compound Poisson
distributions were used. Moreover, as was shown by Arak
(see~\cite{49})the estimate is correct in order. In 1986,  a
joint monograph by Arak and  Zaitsev~\cite{3},
containing a summary of these results, was published in Proceedings of the Steklov Institute of Mathematics. Later
Zaitsev~\cite{43}
showed that a similar estimate holds in the
multidimensional case, and an absolute constant factor is replaced
by~$c(d)$ depending only on the dimension~$d$.

An important
special case of estimating the accuracy of infinitely divisible
approximation is obtained for $\tau=0$, where the right-hand side of
the estimate of  Kolmogorov's uniform distance between
distribution functions
$\rho(\,\cdot\,,\,\cdot\,)$ has the
form $c(d)p$.
In a paper of Zaitsev~\cite{47}, this result
is interpreted as a general estimate for the accuracy of
approximation of the sample composed of non-i.i.d.\ rare events by
a Poisson point process.

In other papers (see~\cite{41} and~\cite{46}), some optimal
bounds for the Kolmogorov distance were also obtained in the
general case. In particular, in the one-dimensional case, they
include simple results which imply simultaneously estimates for the rate of
approximation of convolutions by accompanying infinitely divisible
compound Poisson distributions, and rather general bounds in the
CLT, both optimal in order. Since here the tails of the distributions of the summands are arbitrary, the results cover the now popular case of the so-called
heavy tailed distributions as well.

Similar methods were also used to obtain the following paradoxical
result. There exists a
value $c(d)$} (depending only on the dimension $d$) such that, for any symmetric distribution~$F$
and any~$n\in \mathbf{N}$ the uniform distance between the degrees in the
convolution sense $F^n$ admits the estimates $\rho(F^n,F^{n+1})\le c(d)n^{-1/2}$ and $\rho(F^n,F^{n+2})\le c(d)n^{-1}$, and both estimates are unimprovable
in order (see Zaitsev~\cite{42}).

Now we will apply  Theorem~$\ref{t4}$ and Lemma~$\ref{l1}$ to
obtain the inverse principle type results in the
Littlewood--Offord problem. It is interesting that, in the
multivariate case, the results are obtained by an
application of the one-dimensional  Theorem~\ref{t4} to the
distributions of coordinates of the vector with
distribution~$H_1^{p(1)}$.
\vskip-1mm 

\begin{theorem}
\label{t5}
Let $X$ be a real random variable satisfying
 condition~\eqref{eq16} with~$C_1=1$,
$C_2=\infty$ and $C_3=p(1)>0$. Let $\tau_j\ge\delta_j\ge0$
 and $q_j=Q(F_a^{(j)},\tau_j)$, $j=1,\dots,d$.
Then there exist $r_1,\dots,r_d\in\mathbf N$ and vectors
$u^{(j)}=(u_1^{(j)},\dots, u_{r_j}^{(j)})\in {{\mathbf
R}^{r_j}}$, $j=1,\dots,d$, such that
\begin{equation}
R=\sum_{j=1}^dr_j\ll\sum_{j=1}^d\biggl(|\ln q_j|+\ln\biggl(\frac{\tau_j}{\delta_j}\biggr)+1\biggr)
\label{eq24}
\end{equation}

\vskip-2mm 

\noindent and
\begin{equation}
p(1)M^*\bigl\{\mathbf R^d\setminus\times
_{j=1}^d[K_1(u^{(j)})]_{\delta_j}\bigr\}\ll\sum_{j=1}^d
\biggl(|\ln q_j|+\ln\biggl(\frac{\tau_j}{\delta_j}\biggr)+1\biggr)^3,
\label{eq25}
\end{equation}

\vskip-2mm 

\noindent where $K_1(u^{(j)})\in\mathcal{K}_{r_j}^{(1)}$ and $M^*=\sum_{k=1}^n(E_{a_k}+E_{-a_k})$.

\goodbreak 

Furthermore, the set $\times _{j=1}^dK_1(u^{(j)})$ can be
represented as $K_1(u)\in\mathcal{K}_{R}^{(d)}$,
$u=(u_1,\dots, u_R)\in {(\mathbf R^d)}^R$. Moreover, the vectors $u_s\in\nobreak\mathbf R^d$, 
$s=1,\dots,R$, have  only one non-zero coordinate
each. Denote
$$
s_0=0\quad\text{and}\quad s_k=\sum_{j=1}^{k}r_j,\qquad k=1,\dots,d.
$$
{\bad For $s_{k-1}\,{<}\,s\,{\le}\,s_k$, the vectors $u_s$ are non-zero in the
$k$-th coordinates only and these coordinates are equal to the
sequence of  coordinates $u_1^{(k)},\dots, u_{r_k}^{(k)}$ of the
vector~$u^{(k)}$.} 
\end{theorem}

\begin{theorem}
\label{t6}
Let $X$ be a real random variable satisfying
 condition~\eqref{eq16} with~$C_1=1$, $C_2=\infty$ and
$C_3=p(1)>0$. Let $A,B>0$, $\tau_j\ge \delta_j\ge0$,
$\tau_j/\delta_j\le n^B$ and
 $q_j=Q(F_a^{(j)},\tau_j)\ge n^{-A}$, for $j=1,\dots,d$.
Then there exist numbers  $r_1,\dots,r_d\in\mathbf N$ and vectors
$u^{(j)}=(u_1^{(j)},\dots, u_{r_j}^{(j)})\in {{\mathbf
R}^{r_j}}$, $j=1,\dots,d$, such that
\begin{equation}
R=\sum_{j=1}^dr_j\ll d((A+B)\ln n+1)
\label{eq26}
\end{equation}
and
\begin{equation}
p(1)M^*\bigl\{\mathbf R^d\setminus\times_{j=1}^d[K_1(u^{(j)})]_{\delta_j}\bigr\}\ll d((A+B)\ln n+1)^3,
\label{eq27}
\end{equation}
where $K_1(u^{(j)})\in\mathcal{K}_{r_j}^{(1)}$ and
$M^*=\sum_{k=1}^n(E_{a_k}+E_{-a_k})$. Moreover, the description of the set  $K_1(u)=\times_{j=1}^dK_1(u^{(j)})$ from the end of the formulation of
Theorem~$\ref{t5}$ remains true.
\end{theorem}

\begin{theorem}
\label{t7}
The statements of Theorems\/ $\ref{t5}$ and\/ $\ref{t6}$ remains true with replacing $p(1)$ by $p (0)$
in a particular case, where the
parameters $\tau_j$, $j=1,\dots,d$, involved in the formulations
of these theorems, are all zero.
\end{theorem}

\begin{remark}
\label{r2}\rm
In Theorems
\ref{t5}--\ref{t10}, we use the agreement
 $0/0=1$.
\end{remark}

It is easy to see that, in conditions of Theorem~\ref{t6}
with~$\tau_j=\delta_j\,n^B=\tau$, $j=1,\dots,d$, the set $K_1(u)$
is a GAP of rank  $R=O(\ln n)$, of volume  $3^R=O(n^D)$ (with a
constant~$D$), and such that at least $n-O((\ln
n)^3)$ elements of $a=(a_1,\dots,a_n) \in ({\mathbf{R}^d})^n$
are $\tau/n^B$-close to~$K_1(u)$. Theorem~\ref{t5}
provide
bounds with replacing $\ln n$ by $\left|\ln q\right|$ and without
the assumption~$q=Q(F_a, \tau)\ge n^{-A}$.
Moreover, in~\eqref{eq25} and~\eqref{eq27}, the dependence of constants on
$C_3=p(1)$ is stated explicitly.

Notice that if $\tau_1=\cdots=\tau_d=\tau$, then $q=Q(F_a,
\tau)\le q_j$ and $\left|\ln q_j\right|\le\left|\ln q\right|$,
$j=1,\dots,d$. Moreover, there exist distributions for which the
quantity $q$ may be sufficiently smaller than $\max_j q_j$.
Consider, for instance, the uniform distribution on the boundary
of the square  $\big\{x\in\mathbf R^2\colon |x|=1\big\}$.

In the present paper, we prove as well Theorem~\ref{t8}
which is a multivariate generalization
of Theorem~\ref{t4}.
Furthermore, we state Theorems~\ref{t9} and~\ref{t10}, which are generalizations of Theorems~3.1 and~3.2 of
Chapter~II from~\cite{3}. Deducing Theorems~\ref{t9} and~\ref{t10} from their one-dimensional versions is immediate by
repeating the proof of Theorem~\ref{t8}. Therefore, their
proofs are omitted.

\begin{theorem}
\label{t8}
Let $D$ be a $d$-dimensional
infinitely divisible distribution with characteristic function of
the form  $\exp\{\alpha(\widehat
W(t)\,{-}\,1)\}$, 
$t\,{\in}\,\mathbf R^d$, 
where $\alpha\,{>}\,0$ 
 and~$W$ is a $d$-dimensional probability
distribution. Let
$\tau_j\ge\delta_j\ge0$
 and $\gamma_j=Q(D^{(j)}, \tau_j)$, $j=1,\dots,d$.
Then there exist $r_1,\dots,r_d\in\mathbf N$ and vectors
$u^{(j)}=(u_1^{(j)},\dots, u_{r_j}^{(j)})\in {{\mathbf
R}^{r_j}}$, $j=1,\dots,d$, such that
\begin{equation}
R=\sum_{j=1}^dr_j\ll\sum_{j=1}^d\biggl(|\ln\gamma_j|+\ln\biggl(\frac{\tau_j}{\delta_j}\biggr)+1\biggr)
\label{eq28}
\end{equation}
and
\begin{equation}
\alpha W\bigl\{\mathbf R^d\setminus\times
_{j=1}^d[K_1(u^{(j)})]_{\delta_j}\bigr\}\ll\sum_{j=1}^d
\biggl(|\ln \gamma_j|+\ln\biggl(\frac{\tau_j}{\delta_j}\biggr)+1\biggr)^3,
\label{eq29}
\end{equation}
where $K_1(u^{(j)})\in\mathcal{K}_{r_j}^{(1)}$.
\end{theorem}

\begin{theorem}
\label{t9}
Let  $F_k$, $k=1,\dots,n$, be
$d$-dimensional probability distributions.  Let
$\tau_j\ge\delta_j\ge0$
 and $\gamma_j=Q\bigl(\prod_{k=1}^n
F_k^{(j)}, \tau_j\bigr)$, $j=1,\dots,d$. Then there exist
$r_1,\dots,r_d\in\mathbf N$ and vectors
$u^{(j)}=(u_1^{(j)},\dots, u_{r_j}^{(j)})\in {{\mathbf
R}^{r_j}}$, $j=1,\dots,d$, $x_1,\dots, x_n\in {\mathbf R^d}$,
such that
\begin{equation}
R=\sum_{j=1}^dr_j\ll\sum_{j=1}^d\biggl(|\ln\gamma_j|+\ln\biggl(\frac{\tau_j}{\delta_j}\biggr)+1\biggr)
\label{eq30}
\end{equation}
and
\begin{equation}
\sum_{j=1}^n F_j\bigl\{\mathbf R^d\setminus\times
_{j=1}^d[K_1(u^{(j)})]_{\delta_j}+x_j\bigr\}\ll\sum_{j=1}^d
\biggl(|\ln \gamma_j|+\ln\biggl(\frac{\tau_j}{\delta_j}\biggr)+1\biggr)^3,
\label{eq31}
\end{equation}
where $K_1(u^{(j)})\in\mathcal{K}_{r_j}^{(1)}$.
\end{theorem}

\begin{theorem}
\label{t10}
Let $n\in\mathbf N$ and let $F$ be a
$d$-dimensional probability distribution.  Let
$\tau_j\ge\delta_j\ge0$
 and $\gamma_j=Q((F^{(j)})^n,\tau_j)$, $j=1,\dots,d$.
Then there exist $r_1,\dots,r_d\in\mathbf N$ and vectors
$u^{(j)}=(u_1^{(j)},\dots, u_{r_j}^{(j)})\in {{\mathbf
R}^{r_j}}$, $j=1,\dots,d$, such that
\begin{equation}
R=\sum_{j=1}^dr_j\ll\sum_{j=1}^d\biggl(|\ln\gamma_j|+\ln\biggl(\frac{\tau_j}{\delta_j}\biggr)+1\biggr)
\label{eq32}
\end{equation}
and
\begin{equation}
nF\bigl\{\mathbf R^d\setminus\times
_{j=1}^d[K_1(u^{(j)})]_{\delta_j}\bigr\}\ll\sum_{j=1}^d
\biggl(|\ln \gamma_j|+\ln\biggl(\frac{\tau_j}{\delta_j}\biggr)+1\biggr)^3,
\label{eq33}
\end{equation}
where $K_1(u^{(j)})\in\mathcal{K}_{r_j}^{(1)}$.
\end{theorem}

\begin{remark}
\label{r3}\rm
In Theorems~\ref{t8}, \ref{t9} and~\ref{t10}, the description of the set
$K_1(u)=\times _{j=1}^dK_1(u^{(j)})$ is identical to that
 given at the end of the
formulation of Theorem~$\ref{t5}$.
\end{remark}

\section{Comparison with the results of Nguyen, Tao and
Vu}
\label{s2}

Now we formulate the results discussed in  a review of Nguyen and
Vu \cite{28} (see Theorems~\ref{t11}, \ref{t12} and~\ref{t13}).

A set $K\subset\mathbf R^d$ is called there  \textit{Generalized
Arithmetic Progression} (GAP) of rank $r$ if it can be expressed in
the form
$$
K=\{g_0+m_1g_1+\dots+m_rg_r\colon L_j \le m_j \le L'_j,\,
m_j\in\mathbf Z \text{ for all }1 \le j\le r\},
$$
for some $g_0,\dots, g_r\in\mathbf R^d$, $L_1,\dots, L_r$,
$L'_1,\dots, L'_r\in\mathbf R$.

In fact, $K$ is the image of an integer box
$B=\{(m_1,\dots, m_r) \in\mathbf Z^r\colon  L_j\le m_j \le L'_j\}$
under the linear map
$$
\Phi\colon (m_1,\dots, m_r)\in\mathbf Z^r\to g_0+m_1g_1+\dots+m_rg_r.
$$
The numbers $g_j$ are \textit{generators} of $K$, the numbers $L_j,
L'_j$ are
\textit{dimensions} of~$K$, and   $\operatorname{Vol}(K)=|B|$ is \textit{volume} of~$K$.

We say that $K$ is proper if the map
$\Phi$ is one to one, or,
equivalently, if
$|K|=\operatorname{Vol}(K)$. For non-proper GAPs, we of
course have the strict inequality $|K| < \operatorname{Vol}(K)$. While $-L_j=
L'_j$ for all $j\ge1$ and $g_0=0$, we say that $K$ is symmetric.

First results were related to the discrete case. A few years ago
Tao and Vu~\cite{35} formulated the so-called inverse principle, stating that \textit{a set} $a=(a_1,\dots,a_n)$ \textit{with large small
ball probability must have strong additive structure.}

Here "large small ball probability" means that
$$
Q(F_a,0)=\max_x{\mathbf P}\{S_a=x\}\ge n^{-A}
$$
with some constant $A>0$. "Strong additive structure" means that a large part of
vectors
$a_1,\dots,a_n$ is contained in a GAP with bounded volume. The
following Theorem~\ref{t11} was obtained by Tao and Vu~\cite{35}. In~\cite{27}, this theorem is named "weak
inverse principle" since the choice of $C$ is not optimal.

\begin{theorem}
\label{t11}
Let $0\,{<}\,\varepsilon\,{<}\,1$, 
$A\,{>}\,0$ 
be constants. Then there exist
constants~$r$ and~$C$ depending on~$\varepsilon$ and~$A$ such
that the following holds. Suppose that
$a\,{=}\,(a_1,\dots,\allowbreak a_n) \in ({\mathbf{R}^d})^n$ 
is a multivector in ${\mathbf R}^d$ such
that $Q(F_a, 0)\ge
n^{-A}$. Then there exists a symmetric proper GAP $K$ of constant rank $r$ and of volume  $|K|$ at most $n^C$  such that at least $  n^{1-\varepsilon} $ coordinates of $a$ are
contained in~$K$ $($counting multiplicity$)$.
\end{theorem}

{\bad Later, Tao and Vu~\cite{37} improved the result of
Theorem~\ref{t11}. Nguyen and Vu~\cite{27} have extended the inverse principle to the continuous case (where
$Q(F_a, 0)$} 
 is replaced by $Q(F_a, \tau)$, ${\tau > 0}$) proving,
in particular, the following results.

\begin{theorem}
\label{t12}
Let $X$ be a real random variable satisfying condition~\eqref{eq16} with positive constants $C_1$, $C_2$, $C_3$.
Let $0 < \varepsilon < 1$, $A>0$ be constants and $\tau > 0$ be a parameter that may depend
on~$n$. Suppose that $a=(a_1,\dots,a_n) \in
({\mathbf{R}^d})^n$ is a multivector in~$\mathbf R^d$ such that
 $q=Q(F_a,\tau)\ge n^{-A}$.  Then there exists a symmetric proper GAP $K$ of constant rank $r \ge d$ and of size $|K| =
O(q^{-1}n^{(-r+d)/2})$ such that all but $ \varepsilon n $
coordinates of $a$ are $O( \tau n^{-1/2}\ln
n)$-close to~$K$.
 \end{theorem}

\begin{theorem}
\label{t13}
Let the conditions of Theorem~$\ref{t12}$ be satisfied. Then,
for any number $n'$ between $n^\varepsilon$ and $n$,
 there exists a symmetric proper GAP\enskip
$K=\big\{\sum_{j=1}^rm_jg_j\colon |m_j|\le L_j,\ m_{j}\in\mathbf Z\big\}$ such that

{\rm 1)} At least $n-n'$ elements of $a$ are  $\tau$-close to $K$;

{\rm 2)}~$K$ has small rank $r=O(1)$, and small cardinality
\begin{equation}
|K|\le \max\{O(q^{-1}(n')^{-1/2}), 1\};
\label{eq34}
\end{equation}

{\rm 3)}~There is a non-zero integer $p=O(\sqrt{n'})$ such that all generators $g_j$ of GAP~$K$ have the form $g_j =
(g_{j1},\dots,g_{jd})$, where $g_{jk}=\|a\|\tau p_{jk}/p$
 with~$\|a\|^2=\sum_{j=1}^n\|a_j\|^2$, $p_{jk}\in\mathbf Z$ and $p_{jk}= O(\tau^{-1}\sqrt{n'})$.
\end{theorem}

In the paper~\cite{27}, one can also find some more general
statements, see, for example,~\cite[Theo\-rem~2.9]{27}).

\begin{remark}
\label{r4}\rm
In~\cite{27}, the assumption $\|a\|=1$
 is imposed in the formulations of Theorems~\ref{t12} and~\ref{t13}. Clearly, this assumption can be removed.
\end{remark}

The assertions of Lemma~\ref{l1} and Corollary~\ref{c1} are interesting
only if we assume that $p(\tau/\varkappa)>0$. This condition is
closely related to assumption~\eqref{eq16} in Theorems~\ref{t12}
and~\ref{t13}. Taking into account relations~\eqref{eq2} and
$Q(F_a, \tau)=Q(F_{va}, v\tau)$, $v>0$, we can without loss of generality take
 in~\eqref{eq16} $C_1=1$ and $C_3=p(1)$. Moreover, in our results, $C_2=\infty$. We think that using Lemma~\ref{l1} one could show that $C_2$ may be taken as $C_2=\infty$
 in Theorems~\ref{t12}
and~\ref{t13} too.
Note, however,
 that $p(1)$ is involved in our inequalities explicitly, in contrast with Theorems~\ref{t12} and~\ref{t13}.

Theorem~\ref{t3}
implies Theorem~\ref{t11} and a
one-dimensional version of the first two statements of
Theorem~\ref{t13}.

Thus the following questions arise: what is the relation between
GAPs and CGAPs? Are the assertions about proper GAPs comparable
with the statements concerning CGAPs? In particular, is it
possible to compare Theorems~\ref{t3} and~\ref{t13}? A
positive  answer is given in Proposition~\ref{p1} below.

\begin{proposition}
\label{p1}
Every  one-dimensional CGAP of rank $r$ and volume~$m$ is
contained in a proper symmetric GAP of rank $\le r$ and
volume~$\ll_r m$.
\end{proposition}

Proposition \ref{p1} follows from Theorems~1.6 and 1.9 of Tao and
Vu~\cite{34}. It implies the following Corollary~\ref{c2}.

\begin{corollary}
\label{c2}
Let $K\in\mathcal{K}_{r,m}^{(d)}$ be a $d$-dimensional CGAP of the
form
$K=\times_{j=1}^d K_j$, where $K_j\in\mathcal{K}_{r_j,m_j}$,
$r=(r_1, \dots,r_d)\in\mathbf N^d$, $m=(m_1,
\dots,m_d)\in\mathbf N^d$, with rank  $R=r_1+\dots +r_d$ and volume $M$. Then there exists a proper
$d$-dimensional symmetric GAP $K_0$ of rank $\le R$ and
volume~$\ll_{r,d} M$ and such that
$K\subset K_0$.
\end{corollary}

Thus, inequality~\eqref{eq34} of Theorem~\ref{t13} and
inequality~\eqref{eq19} of Theorem~\ref{t3} (with
$\rho_{n}=1$) are not only of the same form, but their contents
are almost the same, at least for
$d=1$. Attentive readers
may notice evident differences though. In particular, the last
item of Theorem~\ref{t13} is absent in Theorem~\ref{t3}. On
the other hand, in Theorem~\ref{t3}, we take $C_2=\infty$.

One more difference is that, in our theorems, the approximating
set is not proper. However, this leads to a smaller number of
approximating points. Moreover, if $\tau>0$, then it is obvious
that by small perturbations of generators of a non-proper GAP~$K$,
we can construct a proper GAP~$K^*$ with~$[K]_\tau\subset [K^*]_{2\tau}$, the same volume $\operatorname{Vol}(K^*) =
\operatorname{Vol}(K)$ and the same dimensions. The set $[K^*]_{2\tau}$
approximates the
set $a$ not worse than $[K]_\tau$. Note that, according to~\eqref{eq2}, in the conditions of our results there is no essential
difference between $\tau$ and $2\tau$-neighborhoods.

Furthermore, Proposition~\ref{p1} and Corollary~\ref{c2} imply that
we can replace non-proper CGAPs by larger proper GAPs without
essential changes in our formulations.

\begin{remark}
\label{r5}\rm
Using Proposition~\ref{p1}, we can replace CGAPs by symmetric GAPs
of rank $r$ and of volume~${\le m}$, in the definition of
$\mathcal{K}_{r,m}$ and $\beta_{r,m}(W, \tau)$. Then the
assertions of Theorems~\ref{t1} and~\ref{t2} remain valid with
$\le$ replaced by $\ll_r$ in inequalities~\eqref{eq14} and~\eqref{eq17}.
\end{remark}

 It is obvious that the assertions of Lemma~\ref{l1} and Corollary~\ref{c1} may be treated as  statements
 about the measures $G$ and $M^*$. The same may be said about Theorems~\ref{t12} and~\ref{t13}. Moreover,  in the one-dimensional case,  Theorem~\ref{t1} and
Lemma~\ref{l1}
imply precisely the first
 two assertions of Theorem~\ref{t13} (see Theorem~\ref{t3}).

Sometimes, for $d>1$, inequality~\eqref{eq19} (with~$\rho_{n}=1$) may be even stronger than inequality~\eqref{eq34}. For example, if the vector $S_a$ has independent coordinates (this
may happen if each of the vectors~$a_j$ has only one non-zero
coordinate), then
\begin{equation}
q=Q(F_a,\tau)\asymp_d \prod_{j=1}^d q_j.
\label{eq35}
\end{equation}
Note, however, that we could derive a multivariate analogue of
Theorem~\ref{t13} from its one-dimensional version arguing
precisely as in the proof of our Theorem~\ref{t3}.
Then we get
inequality~\eqref{eq19} instead of~\eqref{eq34}.

Theorem~\ref{t3} can be considered as an analogue of both Theorems~\ref{t12} and~\ref{t13}.
Comparing these theorems, we
should mention that the number of approximating points is
sometimes a little bit smaller in Theorem~\ref{t12} than in Theorem~\ref{t3}, but, in
Theorem~\ref{t3}, $C_2=\infty$, and we get a  variety of
results by choosing various $\rho_n$, while in Theorem~\ref{t12}\enskip $\rho_n=n^{-1/2}\ln n$,  and in Theorem~\ref{t13} $\rho_n=1$.

The assertion of Theorem~\ref{t6} implies that, in conditions of
Theorem~\ref{t13}, there exists a symmetric GAP $K$ of rank  $R=O(\ln
n)$,of volume~$3^R=O(n^D)$ and such that at least $n-O((\ln
n)^3)$ elements of  $a=(a_1,\dots,a_n) \in
({\mathbf{R}^d})^n$ are $\tau/n^B$-close to~$K$. Moreover,
Theorem~\ref{t5} provide bounds with replacing $\ln n$ by 
$|\ln q|$ without assumption~$q=Q(F_a, \tau)\ge n^{-A}$
 (recall that this assumption is also absent in conditions of Theorem~\ref{t2}). Comparing with Theorems~\ref{t11},~\ref{t12}
and~\ref{t13}, we see that in Theorem~\ref{t6} the exceptional
set has logarithmic size (which is much better than $O(n)$ and $O(n^\theta)$, $0<\theta\le1$, in Theorems~\ref{t12}
and~\ref{t13}), but this is
attained at the expense of logarithmic growth of the rank.

Notice that all the sets $K_1(u)$ from Theorems~\ref{t5}--\ref{t10}
are simultaneously symmetric GAPs and CGAPs
of rank $R$, and  of volume~$3^R$.

\begin{remark}
\label{r6}\rm
It follows from the proof
that, in Theorem~\ref{t3}, all generators of the GAP
corresponding to $K$ have only one non-zero coordinate each.
\end{remark}

\section{Proof of Theorem~\ref{t3}}
\label{s3}

We will use the classical Ess\'een inequalities  (\cite{14}, see also~\cite{22} and \cite{29}).
\begin{lemma}
\label{l3}
Let $\tau>0$ and
 let $F$ be a $d$-dimensional probability distribution. Then
\begin{equation}
Q(F,\tau)\ll_d \tau^d\int_{|t|\le1/\tau}|\widehat{F}(t)| \,dt,
\label{eq36}
\end{equation}
where $\widehat{F}(t)$ is the characteristic function of the
corresponding random vector.
\end{lemma}

Hal\'asz \cite{21} was the first who has used Ess\'een
inequalities in the Littlewood--Offord problem.

In the general case $Q(F,\tau)$ cannot be estimated from below by the right hand side of
 inequality~\eqref{eq36}. However, if we assume additionally that the distribution $F$ is symmetric
 and its characteristic function is non-negative for all~$t\in\mathbf R$, then we have the lower bound:
\begin{equation}
\label{eq37}
Q(F, \tau)\gg_d \tau^d\int_{|t|\le1/\tau}{\widehat{F}(t)\,dt}
\end{equation}
and, therefore,
\begin{equation}
\label{eq38}
Q(F, \tau)\asymp_d \tau^d\int_{|t|\le1/\tau}{\widehat{F}(t)\,dt}
\end{equation}
(see \cite{1} or \cite [Lemma~1.5 of Chapter II]{3} for $d=1$).
A multidimensional version can be found in \cite{40}, see also \cite{9}.
Using the relation~\eqref{eq38} allowed us to simplify the arguments of Friedland and Sodin~\cite{18},
Rudelson and Vershynin~\cite{32} and Vershynin~\cite{38}, in their studies of the Littlewood--Offord problem
 (see~\cite{9},
\cite{10} and~\cite{12}).

\begin{proof}[Proof of Lemma~\ref{l1}]
Represent the
distribution $G=\mathcal{L}(\widetilde{X})$ as the mixture
$$
G=p_0G_0+p_1G_1,\quad\text{where}\quad p_j={\mathbf P}\bigl\{\widetilde{X} \in A_j\bigr\},\quad j=0,1,
$$
$A_0=\{x\colon |x|\le\tau/\varkappa\}$, $A_1=\{x\colon |x|>\tau/\varkappa\}$,   $G_j$ are probability measures
defined for $p_j>0$ by the formula
$G_j\{B\}=G\{B\cap A_j\}/{p_j}$ , for any Borel set~$B$.
In fact, $G_j$ is the conditional distribution of $\widetilde X$
given that $\widetilde X\in A_j$. If $p_j=0$, then we can take as
$G_j$  an arbitrary  measure. Note that $p_1=p(\tau/\varkappa)$.

For the characteristic function $\widehat{W}(t)=\mathbf{E}\exp(i\langle t,{Y}\rangle)$ of a random vector  $Y\in \mathbf{R}^d$, we
have
$$
|\widehat{W}(t)|^2=\mathbf{E}\exp(i\langle t,\widetilde{Y}\rangle)
= \mathbf{E}\cos(\langle t,\widetilde{Y}\rangle),
$$
where $\widetilde{Y}$ is a corresponding symmetrized random vector. Hence,
\begin{equation}
\label{eq39}
|\widehat{W}(t)| \le
\exp\biggl(-\frac12\bigl(1-|\widehat{W}(t)|^2\bigr)\biggr)=
\exp\biggl(-\frac12\,\mathbf{E}\bigl(1-\cos(\langle t,\widetilde{Y}\rangle)\bigr)\biggr).
\end{equation}

 According to inequalities~\eqref{eq36} and~\eqref{eq39}, we have
\begin{align}
Q(F_a,\tau)&\ll_{d}\tau^d\int_{|t|\le1/\tau}|\widehat{F}_a(t)|\,dt
\nonumber
\\
&\ll_{d}\tau^d\int_{|t|\le1/\tau}\exp\biggl(-\frac12\sum_{k=1}^n\mathbf{E}
\bigl(1-\cos(\langle t,a_k\rangle\widetilde{X})\bigr)\biggr)\,dt=I.
\label{eq40}
\end{align}
It is evident that
\begin{align*}
\sum_{k=1}^n\mathbf{E}\bigl(1-\cos(\langle t,a_k\rangle
\widetilde{X})\bigr)&=\sum_{k=1}^n\int_{-\infty}^{\infty}
\bigl(1-\cos(\langle t,a_k\rangle x)\bigr)G\{dx\}
\\
&=\sum_{k=1}^n\sum_{j=1}^2\int_{A_j}\bigl(1-\cos(\langle t,a_k\rangle x)\bigr)p_jG_j\{dx\}
\\
&\ge\sum_{k=1}^n\int_{A_1}\bigl(1-\cos(\langle t,a_k\rangle x)\bigr)p_1 G_1\{dx\}.
\end{align*}

We now proceed by standard arguments, similarly to the proof of
a result of Ess\'een~\cite{15} (see \cite [Lemma 4 of Chapter II]{29}). Applying Jensen's inequality to the exponential in the integral (see
\cite [p. 49]{29})), we obtain
\begin{align}
I&\le\tau^d\int_{|t|\le1/\tau}\exp\biggl(-\frac{p_1}2\int_{A_1}
\sum_{k=1}^n\bigl(1-\cos(\langle t,a_k\rangle x)\bigr)G_1\{dx\}\biggr)\,dt
\nonumber
\\
&\le\tau^d\int_{|t|\le1/\tau}\int_{A_1}\exp\biggl(-\frac{p_1}2
\sum_{k=1}^n\bigl(1-\cos(\langle t,a_k\rangle x)\bigr)\biggr)G_1\{dx\}\,dt
\nonumber
\\
&\le\sup_{z\in A_1}\tau^d\int_{|t|\le1/\tau}\widehat{H}_z^{p_1}(t)\,dt.
\label{eq41}
\end{align}

Thus, according to~\eqref{eq2} and~\eqref{eq38}, we have
\begin{align}
\sup_{z\in A_1}\tau^d\int_{|t|\le1/\tau}\widehat{H}_z^{p_1}(t)
\,dt &=\sup_{z\ge\tau/\varkappa}\tau^d\int_{|t|\le1/\tau}\widehat{H}_z^{p_1} (t)\,dt \asymp_d\sup_{z\ge\tau/\varkappa}Q(H_z^{p_1},\tau)
\nonumber
\\
&=\sup_{z\ge\tau/\varkappa}Q(H_1^{p_1},\tau/z)=Q(H_1^{p_1},\varkappa),
\label{eq42}
\end{align}
completing the proof.
\end{proof}

\begin{proof}[Proof of Theorem~\ref{t3}]
First we will
 prove Theorem~\ref{t3} for $d=\nobreak 1$. 
Applying  Theorem~\ref{t2} with~$0<\delta=\delta_n=\tau\rho_n\le\tau=\varkappa$ (or with~$\tau=\delta_n=0$, see~\eqref{eq18}), we derive that, for $r,m\in\mathbf N$ the inequality
\begin{equation}
Q(F_a,\tau)\le 2\,c_1^{r+1}\rho_n^{-1}\biggl(\frac1{m\sqrt{\beta_{r,m}(M_0,
\delta_n)}}+\frac{(r+1)^{5r/2}}{(\beta_{r,m}(M_0,\delta_n))^{(r+1)/2}}\biggr) \label{eq43}
\end{equation}
holds, where $M_0=(p(1)/4)M^*$. Let $r=r(A, B,\theta)$
be the minimal
positive  integer such that $A+B<\theta\,(r+1)/2$.
Thus, $r\le\max\{1,2(A+B)/\theta\}$ and
$n^{-A}>n^{B}n^{-\theta(r+1)/2}$ for all~$n>1$. Assume without loss of generality that $n$ is so large that
\begin{align}
n^{-A}&>4c_1^{r+1}n^{B}(r+1)^{5r/2}\biggl(\frac{p(1)\varepsilon n^{\theta}}4\biggr)^{-(r+1)/2}
\nonumber
\\
&\ge4c_1^{r+1}\rho_n^{-1}(r+1)^{5r/2}\biggl(\frac{p(1)\varepsilon n^{\theta}}4\biggr)^{-(r+1)/2}.
\label{eq44}
\end{align}
If~\eqref{eq44} is not satisfied, then $n=O(1)$ and  we can take
as $K$ the set~$K_1(a)\in\mathcal{K}_n^{(1)}$ (see~\eqref{eq20} and~\eqref{eq21}). Choose now a positive integer
$m=\lfloor y\rfloor+1$, where
\begin{equation}
y=\frac{4c_1^{r+1}\rho_n^{-1}}{q\sqrt{p(1)n'/4}}\le m.
\label{eq45}
\end{equation}
Assume that $4p(1)^{-1}\beta_{r,m}(M_0,\delta_n)=\beta_{r,m}(M^*,\delta_n)>n'$.
Recall that $n'\ge\varepsilon n^{\theta}$. Now, using~\eqref{eq43} and our
assumptions, we have
\begin{equation}
n^{-A}\le Q(F_a, \tau) <\frac{Q(F_a, \tau)}2+\frac{n^{-A}}2\le Q(F_a, \tau).
\label{eq46}
\end{equation}
This leads to a contradiction with the assumption $\beta_{r,m}(M^*, \delta_n)>n'$. Hence we conclude that  $\beta_{r,m}(M,\delta_n)\le\beta_{r,m}(M^*,\delta_n)\le n'$.

This means that at least $n-n'$ elements of $a$ are
$\tau\rho_n$-close to a CGAP\enskip $K\in\mathcal{K}_{r,m}$. Equality~\eqref{eq45} implies now relation~\eqref{eq19}.
Theorem~\ref{t3} is proved
 for $d=1$.

Let now $d>1$. We apply Theorem~\ref{t3} with~$d=1$ to the
distributions of the coordinates of the vector~$S_a$,
taking the vector $a^{(j)}=(a_{1j},\dots,a_{nj})$ as vector $a$, for each  $j=1,\dots, d$. Then, for any $a^{(j)}$, there exists a CGAP\enskip $K_j \in \mathcal{K}_{r_j,m_j}$ which satisfies the assertion of Theorem~\ref{t3}, that is:

1)~At least $n-n'$ elements of $a^{(j)}$ are
$\tau\rho_n$-close to~$K_j$;

2)~$K_j$ has small rank $r_j=O(1)$, and
\begin{equation}
\label{eq47}
\begin{gathered}
K_j =\{\langle{\nu}_j,h_j\rangle\colon \nu_j\in\mathbf Z^{r_j}\cap V_j\},\qquad
h_j\in\mathbf R^{r_j},
\\
V_j\subset\mathbf R^{r_j},\quad V_j=-V_j,\quad V_j\text{ is convex},\quad
|\mathbf Z^{r_j}\cap V_j|\le m_j,
\end{gathered}
\end{equation}
where
\begin{equation}
m_j\le\max\bigl\{O\bigl(q_j^{-1}\rho_n^{-1}(n')^{-1/2}\bigr),1\bigr\}.
\label{eq48}
\end{equation}

Thus, the multivector $a^*=(a^{(1)},\dots,a^{(d)})$
is well approximated by the CGAP\enskip $K=\times_{j=1}^d K_j$.
It is easy to see that $K\in\mathcal{K}_{r,m}^{(d)}$,
$r=(r_1,\dots,r_d)\in\mathbf N^d$, $m=(m_1,\dots,m_d)\in{\mathbf
 N}^d$, and
\begin{equation}
\label{eq49}
|\times_{j=1}^d\mathbf Z^{r_j}\cap V|\le\prod_{j=1}^d m_j\le\prod_{j=1}^d
\max\bigl\{O\bigl(q_j^{-1}\rho_n^{-1}(n')^{-1/2}\bigr),1\bigr\},
\end{equation}
where $V=\times _{j=1}^dV_j$.

\goodbreak 

Since at most $n'$ elements of $a^{(j)}$ are far from the CGAPs\enskip $K_j$,
there are at least $n-dn'$ elements of $a$ that are  $\tau\rho_n$-close
to the~\text{CGAP}~$K$. In view of relation~\eqref{eq49} and taking into
account that $K=\times _{j=1}^d K_j$, we obtain
relation~\eqref{eq19}. Theorem~\ref{t3} is proved.
\end{proof}

\begin{remark}
\label{r7}\rm
Notice that, in Theorem~\ref{t3}, the
ranks of $K_j$ are actually the same for all $j=1,\dots,d$. Moreover, in Theorem~\ref{t3}, for sufficiently large~$n$, we
get explicit bound for~$r$: $r\le\max\{1,2(A+B)/\theta\}$.
\end{remark}

\section{Proofs of Theorems~\ref{t5}--\ref{t8}}
\label{s4}

\begin{proof}[Proofs of Theorem~\ref{t5}]
Denote $Q_j=Q(F_a^{(j)},\delta_j)$, $j=1,\dots,d$. By Lemma~\ref{l1} (with~$\varkappa=\tau=\delta_j$),
\begin{equation}
\label{eq50}
Q_j\ll Q(H_{1j}^{p(1)},\delta_j),\qquad j=1,\dots,d,
\end{equation}
where $H_{1j}^{p(1)}$, $j=1,\dots,d$, are the distributions of the coordinates of the vector with distribution $H_1^{p(1)}$. Note that
$H_{1j}^{p(1)}$, $j=1,\dots,d$, are
symmetric infinitely divisible distributions with the L\'evy
 spectral measures~$M_{0j}=(p(1)/4)M^*_j$, where
$M^*_j=\sum_{k=1}^n(E_{a_{kj}}+E_{-a_{kj}})$.

Taking into account~\eqref{eq50}, and applying
Theorem~$\ref{t4}$, we obtain that
 there exist $r_j\in\mathbf N$ and $u^{(j)}=(u_1^{(j)},\dots, u_{r_j}^{(j)})\in {{\mathbf R}^{r_j}}$, $j=1,\dots,d$, such that
\begin{equation}
r_j\ll|\ln Q_j|+1
\label{eq51}
\end{equation}
and
\begin{equation}
p(1)M_j^*\{\mathbf R\setminus[K_1(u^{(j)})]_{\delta_j}\}\ll(|\ln Q_j|+1)^3,
\label{eq52}
\end{equation}
where $K_1(u^{(j)})\in\mathcal{K}_{r_j}^{(1)}$. By~\eqref{eq3},
\begin{equation}
q_j\le(\tau_j/\delta_j+1)Q_j
\label{eq53}
\end{equation}
and
\begin{equation}
|\ln Q_j|\le|\ln q_j|+\ln(\tau_j/\delta_j+1).
\label{eq54}
\end{equation}
Notice that the measures $M_j^*$ are projections of the measure~$M^*$ on the one-dimensional coordinate subspaces.

Constructing now the set $K_1(u)\,{=}\,\times_{j=1}^dK_1(u^{(j)})$ 
as described in the formulation of
Theorem~\ref{t5}, we see that inequalities~\eqref{eq24} and~\eqref{eq25} follow from~\eqref{eq51}, \eqref{eq52} and~\eqref{eq54}. Theorem~\ref{t5} is proved.
\end{proof}

Theorem~\ref{t6} is a direct consequence of Theorem~\ref{t5}.
For the proof of Theorem~\ref{t7} one should replace Lemma~\ref{l1}
by Lemma~\ref{l2} in the proofs of
Theorems \ref{t5} and~\ref{t6}.

\begin{proof}[Proof of Theorem~\ref{t8}]
The proof of Theorem~\ref{t8}  is similar to that of Theorem~\ref{t5}. Recall that the measures
$D^{(j)}$ and $W^{(j)}$, $j=1,\dots,d$, are the projections of
the measures~$D$ and~$W$ respectively on the $j$-th
one-dimensional coordinate subspaces. It is clear that $\widehat D^{(j)}(t)=\exp\{\alpha(\widehat
W^{(j)}(t)-1)\}$, $t\in\mathbf R$. Denote
$\Gamma_j=Q(D^{(j)},\delta_j)$.

Applying Theorem~$\ref{t4}$, we obtain that
 there exist
$r_j\in\mathbf N$ and $u^{(j)}=(u_1^{(j)},\dots,u_{r_j}^{(j)})\in{\mathbf R^{r_j}}$, $j=1,\dots,d$, such that
\begin{equation}
r_j\ll|\ln\Gamma_j|+1
\label{eq55}
\end{equation}
and
\begin{equation}
\alpha W^{(j)}\{{\mathbf R}\setminus[K_1(u^{(j)})]_{\delta_j}\}\ll (|\ln
\Gamma_j|+1)^3,
\label{eq56}
\end{equation}
where $K_1(u^{(j)})\in\mathcal{K}_{r_j}^{(1)}$. By~\eqref{eq3},
\begin{equation}
\gamma_j\le\biggl(\frac{\tau_j}{\delta_j}+1\biggr)\Gamma_j
\label{eq57}
\end{equation}
and
\begin{equation}
|\ln\Gamma_j|\le|\ln\gamma_j|+\ln\biggl(\frac{\tau_j}{\delta_j}+1\biggr).
\label{eq58}
\end{equation}

Defining now the set $K_1(u)\,{=}\,\times_{j=1}^dK_1(u^{(j)})$ 
as described in the formulation of
Theorem~\ref{t5}, we see that inequalities~\eqref{eq28} and~\eqref{eq29} follow from~\eqref{eq55}, \eqref{eq56} and~\eqref{eq58}. Theorem~\ref{t8} is proved.
\end{proof}

\bigskip

We are grateful to a reviewer for useful remarks and for pointing out reference~\cite{34}.


\begin{thebibliography}{20}

\bibitem{1}
 T. V. Arak,   Approximation of $n$-fold convolutions of distributions,
 having a nonnegative characteristic function, with accompanying laws,
Theory Probab. Appl.,  25 (1980), 225--246.

\bibitem{2}
 T. V. Arak,   On the convergence rate in Kolmogorov's uniform limit theorem. I,
Theory Probab. Appl.,  26 (1981), 225--245.

\bibitem{3}
 T. V. Arak, A. Yu. Zaitsev,  Uniform limit theorems for sums of independent random variables,
 Trudy MIAN, 174 (1986), 1--216 (in Russian),
English translation in Proc. Steklov Inst. Math., 174 (1988),
1--216.

\bibitem{4}
S. G. Bobkov, G. P. Chistyakov, Bounds on the maximum of the
density for sums of independent random variables, Zap. Nauchn.
Semin. POMI, 408 (2012), 62--73 (in Russian), English version: J.
Math. Sci. (New York),  199 (2014), 100--106.


\bibitem{5}
V. \v Cekanavi\v cius, Approximation by accompanying distributions
and asymptotic expansions. I,
 Litovsk. Mat. Sb., 29, no. 1 (1989), 171--178 (in Russian).

\bibitem{6}
V.  \v Cekanavi\v cius, Bergstr\"om-type asymptotic expansions in the first uniform Kolmogorov's theorem. In:
 Probability theory and mathematical statistics, Proceedings of the sixth Vilnius conference (Vilnius, 1993),
 Grigelionis, B.  et al. (eds.) Utrecht: VSP, 1994, pp. 223--238.

\bibitem{7}
V.  \v Cekanavi\v cius,  Approximations Methods in Probability
  Theory. Universitext: Springer, 2016, 274 p.

 \bibitem{8}
J.-M. Deshouillers, G. A. Freiman, A. A. Yudin, On bounds for the
concentration function. II, J. Theoret. Probab., 4 (2001),
813--820.


\bibitem{9}
Yu. S. Eliseeva, Multivariate estimates for the concentration
functions of weighted sums of independent identically distributed
random variables, Zap. Nauchn. Semin. POMI, 412 (2013), 121--137
(in Russian), English version: arXiv:1303.4005.


 \bibitem{10}
 Yu. S. Eliseeva, A. Yu. Zaitsev, Estimates for the concentration functions of weighted sums of
 independent random variables, Theory Probab. Appl., 57 (2012), 767--777.

\bibitem{11}
 Yu. S. Eliseeva, A. Yu. Zaitsev, On the Littlewood--Offord problem, Zap.
Nauchn. Semin. POMI, 431 (2014), 72--81 (in Russian), English
version: arXiv:1411.6872 (2014).

\bibitem{12}
Yu. S. Eliseeva, F. G\"otze, A. Yu. Zaitsev, Estimates for the
concentration functions in the Littlewood--Offord problem, Zap.
Nauchn. Semin. POMI, 420 (2013), 50--69 (in Russian), English
version: arXiv:1203.6763 (2012).



\bibitem{13}
P. Erd\"os, On a lemma of Littlewood and Offord,
Bull. Amer. Math. Soc., 51 (1945), 898--902.



\bibitem{14}
 C.-G. Ess\'een, On the Kolmogorov--Rogozin inequality for the
concentration function, Z. Wahrscheinlichkeitstheorie und Verw.
Gebiete,  5 (1966),
 210--216.

\bibitem{15}
C.-G. Ess\'een, On the concentration function of a sum of
independent random variables, Z. Wahrscheinlichkeitstheorie und
Verw. Gebiete, 9 (1968), 290--308.

\bibitem{16}
G. A. Freiman,  Foundations of a structural theory of set addition. Kazan', 1966 (in Russian).
 English translation in: Translations of Mathematical Monographs,
Vol. 37. American Mathematical Society, Providence, R. I., 1973.


\bibitem{17}
G. A. Freiman, A. A. Yudin, On the measure of large values of the
modulus of a trigonometric sum, European J. Combin., 34, no. 8,
(2013), 1338--1347.

\bibitem{18}
O. Friedland, S. Sodin, Bounds on the concentration function
in terms of Diophantine approximation, C. R. Math. Acad. Sci. Paris.
 345 (2007), 513--518.

\bibitem{19}  F. G\"otze,  A. Yu. Zaitsev,
Estimates for the rapid decay of concentration functions of
$n$-fold convolutions, J. Theoret. Probab.  11, no. 3 (1998), 715--731.

\bibitem{20} F. G\"otze,  A. Yu. Zaitsev,
 A multiplicative inequality for concentration functions of $n$-fold convolutions,
 High dimensional probability, v. II (Seattle, WA, 1999),
 Progr. Probab., v. 47, Birkh\"auser Boston, Boston, MA, 2000, pp. 39--47.


\bibitem{21}
G. Hal\'asz, Estimates for the concentration function of combinatorial number theory and probability,
Periodica Mathematica Hungarica, 8 (1977), 197--211.


\bibitem{22}
 W. Hengartner, R. Theodorescu,  Concentration function. Academic Press, New York, 1973.



 \bibitem{23}
A. N. Kolmogorov, Two uniform limit theorems for sums of
independent random variables, Theory Probab. Appl., 1  (1956),
384--394 (in Russian).




\bibitem{24}
J. E. Littlewood, A. C. Offord, On the number of real roots of a
random algebraic equation, Rec. Math. [Mat. Sbornik] N.S., 12
(1943), 277--286.

\bibitem{25}
D. A. Moskvin, G. A. Freiman, A. A. Yudin, Structural theory of
set summation, and local limit theorems for independent lattice
random variables, Theory Probab. Appl., 19 (1974), 52--62 (in
Russian).


\bibitem{26}
 A. B. Mukhin,
The concentration of the distributions of sums of independent
random variables. I; II; III, Izv. Akad. Nauk UzSSR Ser. Fiz.-Mat.
Nauk 17, no. 2 (1973), 25--29; 17, no. 6 (1973), 17--23; 17, no. 1
(1976), 15--19 (in Russian).

 \bibitem{27}
 Hoi Nguyen,  Van Vu, Optimal inverse Littlewood--Offord theorems.
Adv. Math. 226 (2011), 5298--5319.

\bibitem{28}
Hoi Nguyen,  Van Vu, Small ball probabilities, inverse theorems
and applications, Erd\"os Centennial Proceeding, Eds. L. Lov\'asz
et al., Springer, 2013, pp. 409--463, arXiv:1301.0019.


\bibitem{29}
V. V. Petrov,  Sums of independent random variables,  Nauka, Moscow, 1972.

\bibitem{30}
B. A. Rogozin, On the increase of dispersion of sums of
independent random variables, Theory Probab. Appl.,  6 (1961),
106--108.

\bibitem{31}
 M. Rudelson, R. Vershynin,  The Littlewood--Offord problem and
invertibility of random matrices, Adv. Math.,     218 (2008),
600--633.

\bibitem{32}
M. Rudelson, R. Vershynin, Smallest singular value of a random
rectangular matrix, Comm. Pure Appl. Math.,    62 (2009),
1707--1739.

\bibitem{33}
M. Rudelson, R. Vershynin, Small ball probabilities for linear
images of high dimensional distributions,  arXiv:1402.4492 (2014),
to appear in Int. Math. Res. Not.


\bibitem{34}
T. Tao,  Van Vu, John-type theorems for generalized arithmetic
progressions and iterated sumsets, Adv.  Math., 219, no. 2,
(2008), 428--449.

\bibitem{35} T. Tao,  Van Vu, Inverse
Littlewood--Offord theorems and the condition number of random
discrete matrices, Ann. of Math., 169, no. 2 (2009), 595--632.

\bibitem{36}
T. Tao,  Van Vu, From the Littlewood--Offord problem to the
circular law: universality of the spectral distribution of random
matrices, Bull. Amer. Math. Soc. (N.S.),   46 (2009),  377--396.

\bibitem{37}T. Tao,  Van Vu,
A sharp inverse Littlewood--Offord theorem, Random Structures and
Algorithms, 37, no. 4 (2010), 525--539.


\bibitem{38}
R.  Vershynin,
Invertibility of symmetric random matrices,
  Random Structures and Algorithms, 44, no. 2  (2014),  135--182, arXiv:1102.0300.

 \bibitem{39}
 A. Yu. Zaitsev, On the accuracy of approximation of distributions
 of sums of independent random variables -- which are nonzero with a small probability -- by means
 of accompanying laws, Theory Probab. Appl., 28, no. 4 (1984), 657--669.

 \bibitem{40}
 A. Yu. Zaitsev, Multidimensional generalized method of triangular functions,
 Zap. Nauchn. Semin. LOMI, 158 (1987), 81--104 (in Russian).
  English translation in:
  J. Soviet Math., 43, no. 6 (1988), 2797--2810.

 \bibitem{41}
 A. Yu. Zaitsev, On the uniform approximation of distributions of sums of independent random variables,
  Theory Probab. Appl., 32, no. 1 (1987), 40--47.

\bibitem{42}
A. Yu. Zaitsev, Estimates for the closeness of successive convolutions
of multidimensional symmetric distributions,
Probab. Theory  Rel. Fields, 79, no. 2 (1988), 175--200.


\bibitem{43}
A. Yu. Zaitsev, Multivariate version of the second Kolmogorov's uniform limit theorem,
Theory Probab. Appl., 34, no. 1 (1989), 108--128.


 \bibitem{44}
 A. Yu. Zaitsev, On the approximation of convolutions of multi-dimensional
 symmetric distributions by accompaning laws,
 Zap. Nauchn. Semin. LOMI, 177 (1989), 55--72 (in Russian).
  English translation in:
  J. Soviet Math., 61, no. 1 (1992), 1859--1872.

 \bibitem{45}
 A. Yu. Zaitsev, Certain class of nonuniform estimates in multidimensional limit theorems,
 Zap. Nauchn. Semin. LOMI, 184 (1990), 92--105 (in Russian).
  English translation in:
  J. Math. Sci. (N. Y.), 68, no. 4 (1994), 459--468.

\bibitem{46}
A. Yu. Zaitsev, Approximation of convolutions of probability
distributions by infinitely divisible laws under weakened moment
restrictions,
 Zap. Nauchn. Sem. POMI, 194 (1992), 79--90 (in Russian).
  English translation in:
J. Math. Sci. (N. Y.), 75, no. 5 (1995),  1922--1930.

\bibitem{47}
A. Yu. Zaitsev, On approximation of the sample by a Poisson point process,
 Zap. Nauchn. Semin. POMI, 298 (2003), 111--125 (in Russian).
  English translation in:
J. Math. Sci. (N. Y.), 128, no. 1 (2005), 2556--2563.




\bibitem{48}
A. Yu. Zaitsev, Bound for the maximal probability in the Littlewood-Offord problem,
Zap. Nauchn. Sem. POMI, 441 (2015), 204--209 (in Russian).  arXiv:1512.00697.

\bibitem{49}
A. Yu. Zaitsev, T. V. Arak, On the rate of convergence in the
second Kolmogorov's uniform limit theorem, Theory Probab. Appl.,
28, no. 2 (1984), 351--374.










\end{thebibliography}
\end{document}